\documentclass[12pt]{article}
\usepackage{graphics}
\usepackage{amsmath}
\usepackage{amssymb}
\usepackage{graphicx}
\usepackage{makeidx}
\usepackage{mathrsfs}
\usepackage{amsfonts}
\usepackage[mathscr]{eucal}
\usepackage{amsmath}

\DeclareGraphicsExtensions{.eps,.pdf,.jpg,.png}

\def\ds{\displaystyle}
\def\RR{\vbox {\hbox to 8.9pt {I\hskip-2.1pt R\hfil}}}

\def\pni{\par\noindent}
\def\vsh{\smallskip}
\def\vs{\medskip}

\def\vsp{\vsh\pni} %% ie. \smallskip + \par

\newtheorem{theorem}{Theorem}

\newtheorem{corollary}[theorem]{Corollary}

\newtheorem{remark}[theorem]{Remark}

\numberwithin{equation}{section}

\begin{document}
%%%%%%%%%%%%%%%%%%%%%%%%%%%%%%%%%%%%%%%%%%%%%%%%%%%%%%%%%%%%%%%%%%%%%%%%
\vskip 0.75truecm
\font\title=cmbx12 scaled\magstep2
\font\bfs=cmbx12 scaled\magstep1
\font\little=cmr10
\begin{center}
{\title On the Laplace transforms  \\
of derivatives
  of special functions \\
with respect to parameters}
 \\  [0.75truecm]
 Sergei ROGOSIN$^{(1)}$
 Filippo GIRALDI$^{(2), (3)}$
 Francesco MAINARDI$^{(4)}$
\\[0.25truecm]
$^{(1)}$  Department of Economics, Belarusian State University,\\ Nezavisimosti ave 4,  220030 Minsk, Belarus; \\
rogosin@bsu.by; rogosinsv@gmail.com
\\ %[0.25truecm]
$^{(2)}$ Section of Mathematics, International Telematic University Uninettuno,
\\ Corso Vittorio Emanuele II 39, 00186, Rome, Italy;
\\ filippo.giraldi@uninettunouniversity.net; giraldi.filippo@gmail.com\\
%\\ % [0.25truecm]
$^{(3)}$  School of Chemistry and Physics, University of KwaZulu-Natal, \\ Westville Campus, Durban 4000, South Africa
\\ %% [0.25truecm]
$\null^{(4)}$ { Dept. of Physics and Astronomy, University of Bologna, and INFN}
\\ { Via Irnerio 46, I-40126 Bologna, Italy}
\\{ francesco.mainardi@unibo.it; mainardi@bo.infn.it; fracalmo@gmail.com}
\vskip 1truecm %% [0.25truecm]
{\bf Version of
 \today}
 \\
 % [0.25truecm]
 {\bf Paper published in Mathematics (MDPI) 2025 (13), 1980 14pp.}
 \\ {\bf DOI: 103390/math13121980}
\end{center}
%%%%%%%%%%%%
\begin{abstract}
This article is devoted to derivation of the Laplace transforms of the derivatives with respect to parameters of certain special functions, namely, the Mittag-Leffler type, Wright and Le Roy type functions. These formulas show interconnection of these functions and lead to better understanding of their behaviour on the real line. These formulas are represented in the convoluted form and reconstructed in a more suitable form by using Efros theorem.

 \vsp {\bf Keywords}:
{Mittag-Leffler function; Wright function; Le Roy type function; differentiation; Laplace transforms; Efros theorem}

\vsp {\bf 2010 Mathematics Subject Classification (MSC)}:
{44A10, 33B15, 33E12, 33C60}

\end{abstract}

%
% \maketitle

%%%%%%%%%%%%%%%%%%%%%%%%%%%%%
%\begin{document}

%%%%%%%%%%%%%%%%%%%%%%%%%%%%%%%%%%%%%%%%%%
%\setcounter{section}{-1} %% Remove this when starting to work on the template.
%\section{How to Use this Template}

\newpage

\section{Introduction}

An interest to the differential properties of special functions with respect to their parameters is recently growing. Respective formulas become the sources for new classes of special functions as well as differential equations of new types. They are important as for better understanding of the behaviour of special functions as for their applications (see, e.g. \cite{Lan00}). A number of articles are devoted to the differentiation of the Bessel-type functions, see e.g.
\cite{ApeKra85},  \cite{Bry16} and references therein. In particular,
the differential equation for derivatives of the Bessel function w.r.t. parameter was used in \cite{Dun17} for derivation of
an integral representation for $\frac{\partial J_{\nu}(z)}{\partial \nu}$.

Recently Apelblat
\cite{Ape20} and   Apelblat and Mainardi proposed \cite{ApeMai22} an approach for differentiation of the Mittag-Leffler and Wright  functions, respectively,
(see, e.g. \cite{GKMR})
 with respect to parameters.
% It was mostly concerns to the differentiation of the Wright function (see, e.g. \cite{GKMR}).
The obtained formulas are derived by using the formal term-by-term differentiation of the series represented those special  functions.
 So, the derivatives with respect to the parameters are again series
and one can ask where these series do converge.

%In \cite{RoGiMa24} the differentiability with respect to parameters of the general ${\mathcal{HG}}$-function generalizing the Fox ${\mathcal{H}}$-function is studied. These function are introduced in \cite{RogDub23,RogDub24}. They are related to the so-called $I$-integral \cite{In-H87}.

In \cite{RoGiMa24} two approaches were proposed for justification of the formulas of differentiation of the special functions (Mittag-Leffler and Le Roy type) with respect to their parameters. These approaches are based on the  uniform convergence concept and used either series representations or integral representations of the functions mentioned. %In \cite{GiMaRo25}, the  sub-products of \cite{RoGiMa24} are determined and described, namely the differential equations satisfied by the considered special functions are formulated.

This article is devoted to derivation of the Laplace transforms (see, e.g., \cite{DebBha15} of these derivatives with respect to the parameters of certain special functions, namely, the Mittag-Leffler type, Wright and Le Roy type functions (see \cite{GKMR}, \cite{Luc19}, \cite{RogDub23}, \cite{RogDub24}). These formulas show the interconnection of these functions and lead to better understanding of their behavior on the real line. In this connection, we have to mention a survey paper \cite{Ape25} in which some examples are presented justifying importance of the Laplace transforms technique in different branches of applied sciences (cf. also \cite{Ape12}).

Special attention is paid to reconstruction of  these formulas by using Efros theorem \cite{Efr35} in Wlodarski form (see \cite{Wlo52}, \cite{Gra04}).

\newpage
Throughout this paper we assume that all parameters ($\alpha, \beta, \gamma, \alpha_j, \beta_j$) are positive. It follows from \cite{RoGiMa24} that under these conditions the corresponding series converge uniformly with respect to these parameters. The detailed proof is presented in Theorem \ref{Thm2.1}. Thus, we can apply the Laplace transforms to the derivatives point-wise.

%%%%%%%%%%%%%%%%%%%%%%%%%%%%%%%%%%%%%%%%%%
\section{Laplace transforms of the derivatives of the Mittag-Leffler functions}
\label{Laplace_M-L}

We start with the two-parametric Mittag-Leffler function (which is the most useful for applications); see, e.g. \cite[Ch. 4]{GKMR}.
\begin{equation}
\label{2parML1}
E_{\alpha,\beta}(z) = \sum\limits_{k=0}^{\infty} \frac{z^k}{\Gamma(\alpha k + \beta)}.
\end{equation}
It is well known that the function $E_{\alpha,\beta}(z)$ is an entire function of the complex variable $z$ for all ${\mathrm{Re}}\, \alpha > 0, \beta\in {\mathbb C}$. To avoid additional technical details we limit our study to the case  $\alpha > 0,\beta > 0$ only.
The most applicable formula for the Laplace transforms of the function (\ref{2parML1}) is the following one, see, e.g. \cite[(4.9.1)]{GKMR}
\begin{equation}
\label{Lap2parML2}
{\mathcal L}\left(t^{\beta - 1} E_{\alpha,\beta}(\lambda t^{\alpha})\right)(s) = \frac{s^{\alpha-\beta}}{s^{\alpha} - \lambda};\;
\lambda\in {\mathbb R};\; {\mathrm{Re}}\, s > 0;\;\;\; |\lambda s^{-\alpha}| < 1.
\end{equation}
 For the formulas of the derivatives of $t^{\beta - 1} E_{\alpha,\beta}(\lambda t^{\alpha})$ we use the results obtained in \cite{RoGiMa24}.
 In the following theorem we prove the formula for the Laplace transform of the derivative with respect to parameter $\alpha$ of the function $t^{\beta - 1} E_{\alpha,\beta}(\lambda t^{\alpha})$.
 \begin{theorem}
 \label{Thm2.1}
Let $\alpha,\beta > 0,$ $\lambda\in {\mathbb R}$. Then the following formula
 \begin{equation}
\label{2parML3a}
%\begin{split}
\begin{array}{ll}
{\ds {\mathcal L}\left(\frac{\partial}{\partial \alpha} \left[t^{\beta - 1} E_{\alpha,\beta}(\lambda t^{\alpha})\right]\right)(s)}
 &\! =\!
{\ds \int\limits_{0}^{\infty} e^{-st} \sum\limits_{k=0}^{\infty} k \lambda^{k} \left[\frac{t^{\alpha k + \beta - 1} \left(\ln t - \psi(\alpha k + \beta)\right)}{\Gamma(\alpha k + \beta)}\right]}
\\
 &= {\ds - \frac{\lambda s^{\alpha-\beta} \ln\, s}{(s^{\alpha} - \lambda)^2}}
%\end{split}
\end{array}
\end{equation}
is valid for all $|\lambda s^{-\alpha}| < 1$.
Here $\psi$ is the so-called digamma function $\psi(t) = \left(\ln\, \Gamma(t)\right)^{'} = \frac{\Gamma^{'}(t)}{\Gamma(t)}$ (see, e.g., \cite[Ch. 5]{NIST}).
 \end{theorem}

 {\bf Proof.} Let us take arbitrary numbers $a, b$, $0 < a < b < \infty$. Suppose that $\alpha\in [a, b]$.

 First, we prove that the derivative in $\alpha$ of the function $t^{\beta - 1} E_{\alpha,\beta}(\lambda t^{\alpha})$ can be calculated term-wise, i.e.
\begin{equation}
\label{2.3.1}
\frac{\partial}{\partial \alpha} \left[t^{\beta - 1} E_{\alpha,\beta}(\lambda t^{\alpha})\right] = \sum\limits_{k=0}^{\infty} \frac{\partial}{\partial \alpha}
\frac{t^{\alpha k + \beta - 1} \lambda^k}{\Gamma(\alpha k + \beta)} = \sum\limits_{k=0}^{\infty} k \lambda^k \frac{t^{\alpha k + \beta - 1} (\ln t - \psi(\alpha k + \beta))}{\Gamma(\alpha k + \beta)}.
\end{equation}
To see this, it is sufficient to show that the last series is converging uniformly with respect to $\alpha\in [a, b]$.

Let us consider the following series
$$
\ln t \cdot t^{\beta - 1}\sum\limits_{k=0}^{\infty} k \lambda^k \frac{t^{\alpha k}}{\Gamma(\alpha k + \beta)}.
$$
We can find a positive integer number $k_0\in {\mathbb N}$ such that $a k > 4$ for all $k\geq k_0$. Since $\alpha\in [a, b]$, and for integer parts of numbers $[\cdot]$ we have $[a k] \geq a k - 1 > 1, [\beta] > 0$, the the following inequalities can be derived:
$$
\sum\limits_{k=k_0}^{\infty} k \lambda^k \frac{t^{\alpha k}}{\Gamma(\alpha k + \beta)} \leq \sum\limits_{k=k_0}^{\infty} k \lambda^k \frac{t^{b k}}{\Gamma(a k + \beta)} \leq
\sum\limits_{k=k_0}^{\infty} k  \frac{(\lambda t^{b})^{k}}{\Gamma([a k] + [\beta])} \leq
$$
$$
\sum\limits_{k=k_0}^{\infty} \frac{k}{[a k] + [\beta]}  \frac{(\lambda t^{b})^{k}}{\Gamma([a k] + [\beta] - 1)} \leq
\sum\limits_{k=k_0}^{\infty} \frac{k}{a k + \beta - 2}  \frac{(\lambda t^{b})^{k}}{\Gamma(a k + \beta - 3)}.
$$
The series
$$
\sum\limits_{k=k_0}^{\infty} \frac{(\lambda t^{b})^{k}}{\Gamma(a k + \beta - 3)} =
(\lambda t^{b})^{k_0} \sum\limits_{k=k_0}^{\infty} \frac{(\lambda t^{b})^{k - k_0}}{\Gamma(a (k - k_0) + a k_0 + \beta - 3)} =
$$
$$
(\lambda t^{b})^{k_0} \sum\limits_{m=0}^{\infty} \frac{(\lambda t^{b})^{m}}{\Gamma(a m + a k_0 + \beta - 3)} = (\lambda t^{b})^{k_0} E_{a, a k_0 + \beta - 3} (\lambda t^{b})
$$
converges for all $\lambda\in {\mathbb R}$, $t\in {\mathbb R}$ and does not depend on $\alpha$. The sequence $\frac{k}{a k + \beta - 2}$ is bounded and monotonically decreasing. It gives the uniform convergence with respect to $\alpha\in [a, b]$ of the series
$\sum\limits_{k=0}^{\infty} k \lambda^k \frac{t^{\alpha k+\beta -1} \ln t}{\Gamma(\alpha k + \beta)}$.

In order to show the uniform convergence with respect to $\alpha\in [a, b]$ of the series
$\sum\limits_{k=0}^{\infty} k \lambda^k \frac{t^{\alpha k+\beta -1} \psi(\alpha k + \beta)}{\Gamma(\alpha k + \beta)}$ we use the following inequality:
\begin{equation}
\label{psi1}
\ln (x + \frac{1}{2}) \leq \psi(x + 1) \leq \ln (x + e^{-\gamma_{0}}),
\end{equation}
valid for all $x > - 1$   (see \cite{Alz97}, cf. \cite{GuoQi14} for the best possible inequalities of such type). Here $\gamma_{0} := \lim\limits_{n\rightarrow \infty} \left(- \ln n + \sum\limits_{k=1}^{n} \frac{1}{k}\right)$ being the Euler-Mascheroni constant. We also apply Stirling-like asymptotic formula (see, e.g., \cite[Appendix A]{GKMR})
\begin{equation}
\label{gamma1}
\Gamma(\alpha k + \beta) = \sqrt{2\pi} (\alpha k)^{\alpha k + \beta - 1/2} e^{-\alpha k} [1 + o(1)],\;\;\; k \rightarrow \infty.
\end{equation}
Thus
$$ \sum\limits_{k=0}^{\infty} k \lambda^k \frac{t^{\alpha k} \psi(\alpha k + \beta)}{\Gamma(\alpha k + \beta)}$$
$$ \leq
\sum\limits_{k=0}^{\infty} k \lambda^k \frac{t^{b k} \psi(b k + \beta)}{\Gamma(a k + \beta)}
 \leq
\sum\limits_{k=0}^{\infty} \frac{k \ln (bk + \beta - 1 + e^{-\gamma_{0}})}{\sqrt{2\pi} \left(\frac{ak}{e}\right)^{ak} (ak)^{\beta - 1/2}[1 + o(1)]} (\lambda t^b)^k.
$$
Direct calculation show that this power series has an infinite radius of convergence. It gives the uniform convergence with respect to $\alpha\in [a, b]$ of the series
${\ds \sum\limits_{k=0}^{\infty} k \lambda^k \frac{t^{\alpha k+\beta -1} \psi(\alpha k + \beta)}{\Gamma(\alpha k + \beta)}}$.

The ext step is to justify the validity of the formula
\begin{equation}
\label{2.3.2}
{\mathcal L} \left(\sum\limits_{k=0}^{\infty} \frac{\partial}{\partial \alpha} \lambda^{k} \frac{t^{\alpha k + \beta - 1}}{\Gamma(\alpha k + \beta)}\right) =
\sum\limits_{k=0}^{\infty} \int\limits_{0}^{\infty} e^{-st} \frac{\partial}{\partial \alpha} \left(\lambda^{k} \frac{t^{\alpha k + \beta - 1}}{\Gamma(\alpha k + \beta)}\right) dt.
\end{equation}
For this it is sufficient to prove that the series
\begin{equation}
\label{2.3.3}
\sum\limits_{k=0}^{\infty} \int\limits_{0}^{\infty} e^{-st} \frac{\partial}{\partial \alpha} \left(\lambda^{k} \frac{t^{\alpha k + \beta - 1}}{\Gamma(\alpha k + \beta)}\right) dt =
\sum\limits_{k=0}^{\infty} \int\limits_{0}^{\infty} e^{-st} k \lambda^k \frac{t^{\alpha k + \beta - 1} (\ln t - \psi(\alpha k + \beta))}{\Gamma(\alpha k + \beta)}
\end{equation}
converges uniformly with respect to $\alpha\in [a, b]$ in a proper domain of the variable $s$. Let us perform the change of variable $s t = u$ in the following integral
\begin{equation}
\label{2.3.4}
%\begin{array}{ll}
{\ds \int\limits_{0}^{\infty} e^{-st} k \lambda^k \frac{t^{\alpha k + \beta - 1} \ln t}{\Gamma(\alpha k + \beta)}}
 =
{\ds \frac{k (\lambda s^{-\alpha})^k}{s^{\beta}}
 \left[
\int\limits_{0}^{\infty} e^{-u} \frac{u^{\alpha k + \beta - 1} \ln u}{\Gamma(\alpha k + \beta)} d u
  - \int\limits_{0}^{\infty} e^{-u} \frac{u^{\alpha k + \beta - 1} \ln s}{\Gamma(\alpha k + \beta)} d u\right]}.
%\end{array}
\end{equation}
Note that
$\int\limits_{0}^{\infty} e^{-u} u^{\gamma} \ln u \,d u = \left(\int\limits_{0}^{\infty} e^{-u} u^{\gamma} d u\right)^{\prime}_{\gamma}$.
Hence
\begin{equation}
\label{2.3.5}
\begin{array}{ll}
{\ds {\int\limits_{0}^{\infty} e^{-u}
\frac{u^{\alpha k + \beta - 1} \ln u} {\Gamma(\alpha k + \beta)}
\, d u}}
&=
 {\ds \frac{1}{\Gamma(\alpha k + \beta)}
  \left(\int\limits_{0}^{\infty} e^{-u} u^{\alpha k + \beta} \,d u \right)^{\prime}_{\gamma=\alpha k + \beta}}
\\
&=
  {\ds \frac{1}{\Gamma(\alpha k + \beta)} \Gamma^{\prime}(\alpha k + \beta)}
= {\ds \psi(\alpha k + \beta)}.
\end{array}
\end{equation}
For the second term in the right hand side of (\ref{2.3.4}) we have
\begin{equation}
\label{2.3.6}
- \int\limits_{0}^{\infty} e^{-u} \frac{u^{\alpha k + \beta - 1} \ln s}{\Gamma(\alpha k + \beta)} d u = - \ln s \int\limits_{0}^{\infty} e^{-u} \frac{u^{\alpha k + \beta - 1}}{\Gamma(\alpha k + \beta)} d u = - \ln s.
\end{equation}
Substituting (\ref{2.3.5}), (\ref{2.3.6}) into (\ref{2.3.3}) we obtain
\begin{equation}
\label{2.3.7}
\sum\limits_{k=0}^{\infty} \int\limits_{0}^{\infty} e^{-st} \frac{\partial}{\partial \alpha} \left(\lambda^{k} \frac{t^{\alpha k + \beta - 1}}{\Gamma(\alpha k + \beta)}\right) dt =
- \sum\limits_{k=0}^{\infty} \frac{k (\lambda s^{-\alpha})^k}{s^{\beta}} \ln s.
\end{equation}
For the series $\sum\limits_{k=0}^{\infty} k (\lambda s^{-\alpha})^k \leq \sum\limits_{k=0}^{\infty} k (\lambda s^{-a})^k$. The last series does not depend on $\alpha$ and converges in the domain $\left|\lambda s^{-a}\right| < 1$. Hence, we can conclude that the series (\ref{2.3.7}) converges uniformly with respect to $\alpha\in [a, b]$ in the domain $\left|\lambda s^{-\alpha}\right| < 1$. Let us denote $\lambda s^{-\alpha} =: v$ and note that
$$\sum\limits_{k=0}^{\infty} k v^k = v \sum\limits_{k=0}^{\infty} \left(v^k\right)^{\prime}_{v}= \frac{v}{(1 - v)^2}.$$
 Therefore,
$$\sum\limits_{k=0}^{\infty} \int\limits_{0}^{\infty} e^{-st} \frac{\partial}{\partial \alpha} \left(\lambda^{k} \frac{t^{\alpha k + \beta - 1}}{\Gamma(\alpha k + \beta)}\right) dt = - \frac{\ln s \cdot \lambda s^{-\alpha}}{s^{\beta} (1 - \lambda s^{-\alpha})^2} = - \frac{\ln s \cdot \lambda s^{\alpha - \beta}}{(s^{-\alpha} - \lambda )^2}.
$$
It completes the proof. $\Box$
\begin{remark}
\label{Rem2.1}
The obtained result can be formally derived by the direct differentiation of the relation (\ref{Lap2parML2}) with respect to $\alpha$. In the Theorem \ref{Thm2.1} we justify the validity of this formal result.
\end{remark}

Below we present a number of other formulas for the Laplace transform of the derivatives with respect to parameters of different special functions. These formulas are justified in a way similar to that demonstrated in the Theorem \ref{Thm2.1}.
\begin{theorem}
 \label{Thm2.2}
Let $\alpha,\beta > 0,$ $\lambda\in {\mathbb R}$. Then the following formula
 \begin{equation}
\label{2parML3b}
\begin{array}{ll}
{\ds {\mathcal L} \left(\frac{\partial}{\partial \beta} \left[t^{\beta - 1} E_{\alpha,\beta}(\lambda t^{\alpha})\right]\right)(s)}
 &=
{\ds \int\limits_{0}^{\infty} e^{-st} \sum\limits_{k=0}^{\infty}  \lambda^{k} \left[\frac{t^{\alpha k + \beta - 1} \left(\ln t - \psi(\alpha k + \beta)\right)}{\Gamma(\alpha k + \beta)}\right]}
\\
& = {\ds  - \frac{ s^{\alpha-\beta} \ln\, s}{s^{\alpha} - \lambda}}
\end{array}
\end{equation}
is valid for all $|\lambda s^{-\alpha}| < 1$.
\end{theorem}

The direct application of the Laplace transforms to the Mittag-Leffler function leads to a little bit more cumbersome formula
\begin{equation}
\label{Lap2parML3}
{\mathcal L}\left(E_{\alpha,\beta}(t)\right)(s) = \sum\limits_{k=0}^{\infty} \frac{k!}{\Gamma(\alpha k + \beta)} s^{-k-1}.
\end{equation}
It follows from the asymptotic formula (\ref{gamma1}) that the series in the right-hand side of (\ref{Lap2parML3}) converge for all $s\in {\mathbb C}\setminus \{0\}$ whenever $\alpha > 1$, and for all $s\in {\mathbb C},\, |s| > 1,$  whenever $\alpha = 1$. If $0 < \alpha < 1$, then this series diverges everywhere. The corresponding formulas of the Laplace transform of the derivatives of $E_{\alpha,\beta}(t)$ are given in the next therorem.

\begin{theorem}
\label{Thm2.3}
Let $\alpha > 1, \beta\in {\mathbb R}$. Then the following formulas of the Laplace transform of the derivatives of the function $E_{\alpha,\beta}(t)$
\begin{equation}
\label{2parML3a_pure}
{\mathcal L} \left(\frac{\partial}{\partial \alpha} \left[E_{\alpha,\beta}(t)\right]\right)(s)
= - \sum\limits_{k=0}^{\infty} \frac{k!\times k\times \psi(\alpha k + \beta)}{\Gamma(\alpha k + \beta)} s^{-k},
\end{equation}
\begin{equation}
\label{2parML3b_pure}
{\mathcal L} \left(\frac{\partial}{\partial \beta} \left[E_{\alpha,\beta}(t)\right]\right)(s)
= - \sum\limits_{k=0}^{\infty} \frac{k!\times \psi(\alpha k + \beta)}{\Gamma(\alpha k + \beta)} s^{-k}.
\end{equation}
are valid for all $s\in {\mathbb C}\setminus \{0\}$.
\end{theorem}

%%%%%%%%%%%%%%%%%%%%%%%%%%%%%%%%%%%%%%%%%%
\section{Laplace transforms of the derivatives of the Wright function}
\label{Laplace_W}

In this section we provide formulas of the Laplace transforms of the derivatives with respect to parameters of the Wright function (which can be considered for the following real values of parameters $\alpha > - 1, \beta\in {\mathbb R}$)
\begin{equation}
\label{Wright1}
\phi(\alpha,\beta;z) = W_{\alpha,\beta}(z) = \sum\limits_{k=0}^{\infty} \frac{z^k}{k! \Gamma(\alpha k + \beta)}.
\end{equation}
Two formulas of the Laplace transforms of this function are known.The first one represents the Laplace transforms in terms of the Mittag-Leffler function (\ref{2parML1}) in the case $\alpha > 0$
\begin{equation}
\label{LapWright1}
{\mathcal L} \left[W_{\alpha,\beta}(\pm t)\right](s) = \frac{1}{s} E_{\alpha,\beta}(\pm s^{-1}).
\end{equation}
As in Theorem \ref{Thm2.1} we can justify the following formulas of the Laplace transforms of the derivatives with respect to parameters:
\begin{theorem}
\label{Thm3.1}
Let $\alpha > 0, \beta > 0$. Then the following formulas for the Laplace transforms of the derivatives of $W_{\alpha,\beta}(t)$ have the form.
\begin{equation}
\label{LapWright1a}
{\mathcal L} \left(\frac{\partial}{\partial \alpha} \left[W_{\alpha,\beta}(\pm t)\right]\right)(s) = - \frac{1}{s} \sum\limits_{k=0}^{\infty} \frac{\pm s^{- k} k \psi(\alpha k + \beta)}{\Gamma(\alpha k + \beta)},
\end{equation}
\begin{equation}
\label{LapWright1b}
{\mathcal L} \left(\frac{\partial}{\partial \beta} \left[W_{\alpha,\beta}(\pm t)\right]\right)(s) = - \frac{1}{s} \sum\limits_{k=0}^{\infty} \frac{\pm s^{- k} \psi(\alpha k + \beta)}{\Gamma(\alpha k + \beta)}.
\end{equation}
are valid for all $s\in {\mathbb C}\setminus \{0\}$.
\end{theorem}

Second variant of the Laplace transform of the Wright function is similar to formulas (\ref{2parML3a}), (\ref{2parML3b}).
\begin{equation}
\label{LapWright2}
{\mathcal L} \left(t^{\beta - 1} \left[W_{\alpha,\beta}(\lambda t^{\alpha})\right]\right)(s) = \frac{1}{s^{\beta}} \sum\limits_{k=0}^{\infty} \frac{(\lambda s^{-\alpha})^k}{k!} = \frac{1}{s^{\beta}} e^{\lambda s^{-\alpha}}.
\end{equation}
\begin{corollary}
\label{Cor3.1}
Let $\alpha > 0, \beta > 0$. Then the following formulas for the Laplace transforms of the derivatives have the form.
\begin{equation}
\label{LapWright2a}
{\mathcal L} \left(\frac{\partial}{\partial \alpha} t^{\beta - 1} \left[W_{\alpha,\beta}(\lambda t^{\alpha})\right]\right)(s) = \frac{1}{s^{\beta}} \sum\limits_{k=0}^{\infty} \frac{(\lambda s^{-\alpha})^k}{k!} = \frac{-\lambda \ln\, s}{s^{\alpha + \beta}} e^{\lambda s^{-\alpha}},
\end{equation}
\begin{equation}
\label{LapWright2b}
{\mathcal L} \left(\frac{\partial}{\partial \beta} t^{\beta - 1} \left[W_{\alpha,\beta}(\lambda t^{\alpha})\right]\right)(s) = \frac{1}{s^{\beta}} \sum\limits_{k=0}^{\infty} \frac{(\lambda s^{-\alpha})^k}{k!} = \frac{- \ln\, s}{s^{\beta}} e^{\lambda s^{-\alpha}}.
\end{equation}
are valid for all $s\in {\mathbb C}, {\mathrm{Re}}\, s > 0$.
\end{corollary}

%%%%%%%%%%%%%%%%%%%%%%%%%%%%%%%%%%%%%%%%%%
\section{Laplace transforms of the derivatives of the three- and four-parametric Mittag-Leffler functions}
\label{Laplace_3-4M-L}

The most useful for applications among the multi-parametric functions is the so called Prabhakar function (see, e.g., \cite[Sec. 5.1]{GKMR})
\begin{equation}
\label{3parML1}
E_{\alpha,\beta}^{\gamma}(z) = \sum\limits_{k=0}^{\infty} \frac{(\gamma)_k z^k}{\Gamma(\alpha k + \beta)} =
\sum\limits_{k=0}^{\infty} \frac{\Gamma(\gamma) z^k}{\Gamma(\gamma + k)\Gamma(\alpha k + \beta)}.
\end{equation}
The Laplace transforms of the Prabhakar function is given by the following formula
\begin{equation}
\label{Lap3parML1}
{\mathcal L}\left(t^{\beta - 1} E_{\alpha,\beta}^{\gamma}(\lambda t^{\alpha})\right)(s) = \frac{s^{\alpha \gamma - \beta}}{(s^{\alpha} - \lambda)^{\gamma}};\;
\lambda\in {\mathbb R};\; {\mathrm{Re}}\, s > 0;\;\;\; |\lambda s^{-\alpha}| < 1.
\end{equation}
In a way similar to that of Theorem \ref{Thm2.1} we can justify the following formulas of the Laplace transforms of the derivatives of this function with respect to parameters.
\begin{theorem}
\label{Thm4.1}
Let $\alpha, \beta, \gamma > 0$. Then the following formulas of the Laplace transform of the derivatives with respect to parameters are satisfied:
\begin{equation}
\label{Lap3parML2a}
{\mathcal L} \left(\frac{\partial}{\partial \alpha} \left[t^{\beta - 1} E_{\alpha,\beta}^{\gamma}(\lambda t^{\alpha})\right]\right)(s) = \frac{- \gamma \lambda s^{\alpha \gamma - \beta} \ln\,s}{(s^{\alpha} - \lambda)^{\gamma + 1}};\;
\lambda\in {\mathbb R};\; {\mathrm{Re}}\, s > 0;\;\;\; |\lambda s^{-\alpha}| < 1.
\end{equation}
\begin{equation}
\label{Lap3parML2b}
{\mathcal L} \left(\frac{\partial}{\partial \beta} \left[t^{\beta - 1} E_{\alpha,\beta}^{\gamma}(\lambda t^{\alpha})\right]\right)(s) = \frac{-  s^{\alpha \gamma - \beta} \ln\,s}{(s^{\alpha} - \lambda)^{\gamma}};\;
\lambda\in {\mathbb R};\; {\mathrm{Re}}\, s > 0;\;\;\; |\lambda s^{-\alpha}| < 1.
\end{equation}
\begin{equation}
\label{Lap3parML2c}
{\mathcal L} \left(\frac{\partial}{\partial \gamma} \left[t^{\beta - 1} E_{\alpha,\beta}^{\gamma}(\lambda t^{\alpha})\right]\right)(s) = \frac{s^{\alpha \gamma - \beta}}{(s^{\alpha} - \lambda)^{\gamma}}\left[- \alpha  \ln\,s +  \ln\,(s^{\alpha} - \lambda)\right];
\end{equation}
$$
\lambda\in {\mathbb R};\; {\mathrm{Re}}\, s > 0;\;\;\; |\lambda s^{-\alpha}| < 1.
$$
\end{theorem}

The second variant of the Laplace has the following form
\begin{equation}
\label{Lap3parML3}
{\mathcal L}\left(t^{\gamma - 1} E_{\alpha,\beta}^{\gamma}(t)\right)(s) = \Gamma(\gamma) \sum\limits_{k=0}^{\infty} \frac{s^{-\gamma - k}}{\Gamma(\alpha k + \beta)} = \frac{\Gamma(\gamma)}{s^{\gamma}} E_{\alpha,\beta}\left(\frac{1}{s}\right);\;
{\mathrm{Re}}\, s > 0.
\end{equation}
\begin{corollary}
Let $\alpha, \beta, \gamma > 0$. Then the
 following formulas of the Laplace transforms of the derivatives with respect to parameters of the function $t^{\gamma - 1} E_{\alpha,\beta}^{\gamma}(t)$ are satisfied for ${\mathrm{Re}}\, s > 0$
\begin{equation}
\label{Lap3parML4a}
{\mathcal L}\left(\frac{\partial}{\partial \alpha} \left[t^{\gamma - 1} E_{\alpha,\beta}^{\gamma}(t)\right]\right)(s) = -  \frac{\Gamma(\gamma)}{s^{\gamma}} \sum\limits_{k=0}^{\infty} \frac{k \psi(\alpha k + \beta)}{\Gamma(\alpha k + \beta)} \frac{1}{s^k}.
\end{equation}
\begin{equation}
\label{Lap3parML4b}
{\mathcal L}\left(\frac{\partial}{\partial \beta} \left[t^{\gamma - 1} E_{\alpha,\beta}^{\gamma}(t)\right]\right)(s) = -  \frac{\Gamma(\gamma)}{s^{\gamma}} \sum\limits_{k=0}^{\infty} \frac{\psi(\alpha k + \beta)}{\Gamma(\alpha k + \beta)} \frac{1}{s^k}.
\end{equation}
\begin{equation}
\label{Lap3parML4c}
{\mathcal L}\left(\frac{\partial}{\partial \gamma} \left[t^{\gamma - 1} E_{\alpha,\beta}^{\gamma}(t)\right]\right)(s) = \frac{\Gamma(\gamma) [\psi(\gamma) - \ln\, s]}{s^{\gamma}} E_{\alpha,\beta}\left(\frac{1}{s}\right).
\end{equation}
\end{corollary}

The results for the four-parametric Mittag-Leffler function
\begin{equation}
\label{4parML1}
E_{(\alpha_1,\beta_1),(\alpha_2,\beta_2)}(z) = \sum\limits_{k=0}^{\infty} \frac{z^k}{\Gamma(\alpha_1 k + \beta_1) \Gamma(\alpha_2 k + \beta_2)}.
\end{equation}
are similar to the obtained above. For instance we can use the following formula for the Laplace transforms ($\lambda\in {\mathbb R};\; {\mathrm{Re}}\, s > 0$):
\begin{equation}
\label{Lap4parML2}
{\mathcal L} \left(t^{\beta_1 - 1} E_{(\alpha_1,\beta_1),(\alpha_2,\beta_2)}(\lambda t^{\alpha_1})\right)(s) = \frac{1}{s^{\beta_1}} \sum\limits_{k=0}^{\infty} \frac{(\lambda s^{-\alpha_1})^k}{\Gamma(\alpha_2 k + \beta_2)} = \frac{1}{s^{\beta_1}} E_{\alpha_2,\beta_2}(\lambda s^{-\alpha_1}).
\end{equation}
\begin{theorem}
\label{Thm4.2}
Let $\alpha_1, \alpha_2, \beta_1, \beta_2 > 0$
Then the following formulas of the Laplace transforms of the derivatives are satisfied for ${\mathrm{Re}}\, s > 0$.
\begin{equation}
\label{Lap4parML3a1}
{\mathcal L} \left(\frac{\partial}{\partial \alpha_1} \left[t^{\beta_1 - 1} E_{(\alpha_1,\beta_1),(\alpha_2,\beta_2)}(\lambda t^{\alpha_1})\right]\right)(s) =
- \frac{1}{s^{\beta_1}} \sum\limits_{k=0}^{\infty} \frac{k (\lambda s^{-\alpha_1})^k}{\Gamma(\alpha_2 k + \beta_2)}.
\end{equation}
\begin{equation}
\label{Lap4parML3b1}
{\mathcal L} \left(\frac{\partial}{\partial \beta_1} \left[t^{\beta_1 - 1} E_{(\alpha_1,\beta_1),(\alpha_2,\beta_2)}(\lambda t^{\alpha_1})\right]\right)(s) =
- \frac{\ln\, s}{s^{\beta_1}} E_{\alpha_2,\beta_2}(\lambda s^{-\alpha_1}).
\end{equation}
\begin{equation}
\label{Lap4parML3a2}
{\mathcal L} \left(\frac{\partial}{\partial \alpha_2} \left[t^{\beta_1 - 1} E_{(\alpha_1,\beta_1),(\alpha_2,\beta_2)}(\lambda t^{\alpha_1})\right]\right)(s) =
- \frac{1}{s^{\beta_1}} \sum\limits_{k=0}^{\infty} \frac{k \psi(\alpha_2 k + \beta_2)}{\Gamma(\alpha_2 k + \beta_2)} (\lambda s^{-\alpha_1})^k.
\end{equation}
\begin{equation}
\label{Lap4parML3b2}
{\mathcal L} \left(\frac{\partial}{\partial \beta_2} \left[t^{\beta_1 - 1} E_{(\alpha_1,\beta_1),(\alpha_2,\beta_2)}(\lambda t^{\alpha_1})\right]\right)(s) =
- \frac{1}{s^{\beta_1}} \sum\limits_{k=0}^{\infty} \frac{\psi(\alpha_2 k + \beta_2)}{\Gamma(\alpha_2 k + \beta_2)} (\lambda s^{-\alpha_1})^k.
\end{equation}
\end{theorem}

%%%%%%%%%%%%%%%%%%%%%%%%%%%%%%%%%%%%%%%%%%
\section{Laplace transforms of the derivatives of the Le Roy type function}
\label{Laplace_Le Roy}

Here we deal with the Laplace transforms of the derivatives of the Le Roy type function (see, e.g., \cite[Sec. 5.3]{GKMR}, \cite{RogDub23}, \cite{RogDub24})
\begin{equation}
\label{LeRoy1}
F_{\alpha,\beta}^{(\gamma)}(z) = \sum\limits_{k=0}^{\infty} \frac{z^k}{[\Gamma(\alpha k + \beta)]^{\gamma}}
\end{equation}
which generalizes the function
\begin{equation}
\label{LeRoy0}
R_{\gamma}(z) = \sum\limits_{k=0}^{\infty} \frac{z^k}{[k!]^{\gamma}},
\end{equation}
introduced by E. Le Roy at the end of XIX century for the study of the analytic continuation of power series.

Let us calculate the Laplace transforms of the Le Roy type function taken in the form similar to (\ref{Lap2parML2}), (\ref{LapWright1a}), (\ref{Lap3parML3})
\begin{equation}
\label{LapLeRoy1}
{\mathcal L}\left(t^{\gamma - 1} F_{\alpha,\beta}^{(\gamma)}(\lambda t^{\alpha})\right)(s) = \frac{1}{s^{\beta}} \sum\limits_{k=0}^{\infty} \frac{(\lambda s^{-\alpha})^k}{[\Gamma(\alpha k + \beta)]^{\gamma - 1}} = \frac{1}{s^{\beta}} F_{\alpha,\beta}^{(\gamma)}\left(\frac{\lambda}{s^{\alpha}}\right);
\end{equation}
$$
\lambda\in {\mathbb R},\; {\mathrm{Re}}\, s > 0;\; \gamma > 1.
$$
In particular, for $\gamma = 2$ the Laplace transforms image is related to the Mittag-Leffler function
\begin{equation}
\label{LapLeRoy2}
{\mathcal L}\left(t^{\gamma - 1} F_{\alpha,\beta}^{(2)}(\lambda t^{\alpha})\right)(s) = \frac{1}{s^{\beta}} E_{\alpha,\beta}\left(\frac{\lambda}{s^{\alpha}}\right).
\end{equation}
\begin{theorem}
\label{Thm5.1} Let $\alpha > 0, \beta > 0, \gamma > 1$.
Then the formulas of the Laplace transforms of the derivatives have the following forms ($\lambda\in {\mathbb R},\; {\mathrm{Re}}\, s > 0$):
\begin{equation}
\label{LapLeRoy3a}
{\mathcal L} \left(\frac{\partial}{\partial \alpha} \left[t^{\gamma - 1} F_{\alpha,\beta}^{(\gamma)}(\lambda t^{\alpha})\right]\right)(s) =
- \frac{1}{s^{\beta}} \sum\limits_{k=0}^{\infty} k \left(\frac{\lambda}{s^{\alpha}}\right)^k \frac{\ln\, s + (\gamma - 1) \psi(\alpha k + \beta)}{[\Gamma(\alpha k + \beta)]^{\gamma - 1}};
\end{equation}
\begin{equation}
\label{LapLeRoy3b}
{\mathcal L} \left(\frac{\partial}{\partial \beta} \left[t^{\gamma - 1} F_{\alpha,\beta}^{(\gamma)}(\lambda t^{\alpha})\right]\right)(s) =
- \frac{1}{s^{\beta}} \sum\limits_{k=0}^{\infty} \left(\frac{\lambda}{s^{\alpha}}\right)^k \frac{\ln\, s + (\gamma - 1) \psi(\alpha k + \beta)}{[\Gamma(\alpha k + \beta)]^{\gamma - 1}};
\end{equation}
\begin{equation}
\label{LapLeRoy3c}
{\mathcal L} \left(\frac{\partial}{\partial \gamma} \left[t^{\gamma - 1} F_{\alpha,\beta}^{(\gamma)}(\lambda t^{\alpha})\right]\right)(s) =
- \frac{1}{s^{\beta}} \sum\limits_{k=0}^{\infty} \left(\frac{\lambda}{s^{\alpha}}\right)^k \frac{\ln\, \Gamma(\alpha k + \beta)}{[\Gamma(\alpha k + \beta)]^{\gamma - 1}}.
\end{equation}
\end{theorem}

%%%%%%%%%%%%%%%%%%%%%%%%%%%%%%%%%%%%%%
\section{Convoluted forms for the derivatives of
Mittag-Leffler and Wright type functions with respect to parameters}
\label{6}

The above-reported Laplace transforms of the partial derivatives of Mittag-Leffler and Wright-type functions with respect to parameters allow us to find new representations of the partial derivatives of these special functions in terms of convoluted forms. These forms are obtained via the convolution theorem and the Efros theorem, which are reported below for the sake of clarity.

As far as the convolution theorem is concerned, let $\varphi_1$ and $\varphi_2$ be functions, locally integrable on $\left[\right. 0, +\infty\left.\right)$, that exhibit exponential growth $u_1$ and $u_2$, respectively. Thus, positive parameters $m_j$, $u_j$, $T_j$, exist such that $\left|\varphi_j\left(t\right)\right|\leq m_j \exp\left(u_j t\right)$, for all $t \geq T_j$, and $j=1,2$. The convolution product $\left(\varphi_1\ast \varphi_2\right)(t)$, given by
\begin{eqnarray}
&&\hspace{-2em}
\left(\varphi_1\ast \varphi_2\right)(t)
=\int_0^t\varphi_1\left(\tau\right)
\varphi_2\left(t-\tau\right)d \tau,
\label{convprod}
\end{eqnarray}
is properly defined on $\left[\right. 0, +\infty\left.\right)$ and exhibits exponential growth $\max \left\{u_1,u_2\right\}$. Additionally, the following Laplace transform:
\begin{eqnarray}
&&\hspace{-2em}
\mathcal{L}\left(\left(\varphi_1\ast \varphi_2\right)(t)\right)(s)
=\mathcal{L}\left(\varphi_1(t)\right)(s)\mathcal{L}\left(\varphi_2(t)\right)(s),
\label{Lconvprod}
\end{eqnarray}
 holds for $\operatorname{Re}\, s > \max \left\{u_1,u_2\right\}$.

At this stage, we are equipped to enunciate the various theorems concerning the representations of the special functions under study mentioned above.

\begin{theorem}
\label{Th1}
The partial derivative $\partial/\partial \alpha \left(
t^{\beta-1} E_{\alpha,\beta}\left(\lambda t^{\alpha}\right)\right)$,
involving the Mittag-Leffler function, results to be a convoluted form with kernel of logarithmic type:
\begin{eqnarray}
&&\hspace{-2em}
\frac{\partial}{\partial \alpha}\,
t^{\beta-1} E_{\alpha,\beta}\left(\lambda t^{\alpha}\right)=\lambda \int_0^t \left(\ln \left(t-t_1\right)+\gamma_0\right) t_1^{\alpha+\beta-2}E^2_{\alpha,\alpha+\beta-1}\left(\lambda t_1^{\alpha}\right)dt_1, \nonumber \\ &&\label{Datbm1Eab}
\end{eqnarray}
for every $t\geq 0$, where $\gamma_0$ is the Euler-Mascheroni constant.
\end{theorem}

For $\alpha>0$ and $\left(\alpha+\beta\right)>1$, the first term of the right side of Eq. (\ref{Datbm1Eab}) is given by the expression below,
\begin{eqnarray}
&&\hspace{-1em}\int_0^t \ln\left(t-t_1\right) t_1^{\alpha+\beta-2}E^2_{\alpha,\alpha+\beta-1}\left(\lambda t_1^{\alpha}\right)dt_1\nonumber \\ &&\hspace{-1em}=
t^{\alpha+\beta-1}\left(\ln t-\gamma_0\right)
\left(\left(\alpha+\beta-1\right)E^2_{\alpha,\alpha+\beta-1}\left(\lambda t^{\alpha}\right)
+2\alpha \lambda t^{ \alpha}E^3_{\alpha,2\alpha+\beta-1}\left(\lambda t^{\alpha}\right)\right)\nonumber \\ &&\hspace{-1em}
-t^{\alpha+\beta-1}
\left(
\left(\alpha+\beta-1\right)
\sum_{n=0}^{\infty}\frac{\left(n+1\right)
\psi\left(\alpha n+\alpha+\beta\right)\left(\lambda t^{\alpha}\right)^n}{\Gamma\left(\alpha n+\alpha+\beta-1\right)}\right.\nonumber \\ &&\hspace{-1em}
\left.+\alpha \lambda t^{\alpha}
\sum_{n=0}^{\infty}\frac{\left(n+2\right)\left(n+1\right)
\psi\left(\alpha n+2\alpha+\beta\right)\left(\lambda t^{\alpha}\right)^n}{\Gamma\left(\alpha n+2\alpha+\beta-1\right)}\right)
.\label{Datbm1Eabrhs1}
\end{eqnarray}

For $\alpha>0$ and $\left(\alpha+\beta\right)>1$, the second term on the right side of Eq. (\ref{Datbm1Eab}) is given by the form below,
\begin{eqnarray}
&&\hspace{-2em}\int_0^t  t_1^{\alpha+\beta-2}E^2_{\alpha,\alpha+\beta-1}\left(\lambda t_1^{\alpha}\right)dt_1=t^{\alpha+\beta-1}E^2_{\alpha,\alpha+\beta}\left(\lambda t^{\alpha}\right).\label{Datbm1Eabrhs2}
\end{eqnarray}

{\bf Proof.} [Proof of Theorem \ref{Th1}]
Form (\ref{Datbm1Eab}) is obtained by expressing the right side of Eq. (\ref{2parML3a}) as follows:
$$ \lambda \frac{- \ln s}{s} \frac{s^{1+\alpha-\beta}}{\left(s^{\alpha}-\lambda\right)^2}.$$
The Laplace inversion of the above-reported expression is obtained by performing the Laplace inversion of each of the two fractions and adopting the convolution theorem \cite{PrBrMa145}.

The term-by-term integration of the power series representation of the involved Mittag-Leffler function provides Eq. (\ref{Datbm1Eabrhs1}) via the following integral \cite{PrBrMa145}:
\begin{eqnarray}
&&\hspace{-1em}\int_0^t t_1^{r} \ln\left(t-t_1\right)  dt_1
=\left(1+r\right)t^{1+r}\left( \ln t -\gamma_0 -\psi\left(2+r\right)\right)
,\label{intplog}
\end{eqnarray}
holding for every $t\geq 0$ and $r>-1$. Instead, Eq. (\ref{Datbm1Eabrhs2}) is obtained in a straightforward way via the term-by-term integration of the power series representation of the Mittag-Leffler function involved.
$\Box$

\begin{theorem}
\label{Th2}
The partial derivative $\partial/\partial \beta \left(
t^{\beta-1} E_{\alpha,\beta}\left(\lambda t^{\alpha}\right)\right)$,
involving the Mittag-Leffler function, results to be a convoluted form with kernel of logarithmic type:
\begin{eqnarray}
&&\hspace{-2em}\frac{\partial}{\partial \beta}
\left(t^{\beta-1} E_{\alpha,\beta}\left(\lambda t^{\alpha}\right)\right)= \int_0^t \left(\ln \left(t-t_1\right)+\gamma_0\right) t_1^{\beta-2}E_{\alpha,\beta-1}\left(\lambda t_1^{\alpha}\right)dt_1, \nonumber \\ &&\label {Dbtbm1Eab}
\end{eqnarray}
for every $t\geq 0$.
\end{theorem}

For $\alpha>0$ and $\beta >1$, the first term on the right side of Eq. (\ref{Dbtbm1Eab}) is given by
\begin{eqnarray}
&&\hspace{-1em}\int_0^t \left(\ln\left(t-t_1\right)\right) t_1^{\beta-2}E_{\alpha,\beta-1}\left(\lambda t_1^{\alpha}\right)dt_1\nonumber \\ &&\hspace{-1em}=
t^{\beta-1}\left(\ln t-\gamma_0\right)
\left(\left(\beta-1\right)E_{\alpha,\beta-1}\left(\lambda t^{\alpha}\right)
+\alpha \lambda t^{ \alpha}E^2_{\alpha,\alpha+\beta-1}\left(\lambda t^{\alpha}\right)\right)\nonumber \\ &&\hspace{-1em}
-t^{\beta-1}
\left(
\left(\beta-1\right)
\sum_{n=0}^{\infty}\frac{
\psi\left(\alpha n+\beta\right)\left(\lambda t^{\alpha}\right)^n}{\Gamma\left(\alpha n+\beta-1\right)}\right.\nonumber \\ &&\hspace{-1em}
\left.+ \alpha
\sum_{n=1}^{\infty}\frac{n
\psi\left(\alpha n+\beta\right)\left(\lambda t^{\alpha}\right)^n}{\Gamma\left(\alpha n+\beta-1\right)}\right)
.\label{Dbtbm1Eabrhs1}
\end{eqnarray}

For $\alpha>0$ and $\beta>1$, the second term of the right side of Eq. (\ref{Dbtbm1Eab}) is given by
\begin{eqnarray}
&&\hspace{-2em}\int_0^t  t_1^{\beta-2}E_{\alpha,\beta-1}\left(\lambda t_1^{\alpha}\right)dt_1=t^{\beta-1}E_{\alpha,\beta}\left(\lambda t^{\alpha}\right).\label{Dbtbm1Eabrhs2}
\end{eqnarray}
The proof of Theorem \ref{Th2} is similar to the proof of Theorem \ref{Th1}.

The Efros theorem is a generalization of the convolution theorem that represents in the present scenario a valuable tool for evaluating the involved Laplace inversions \cite{Efr35}.  Let $F(s)$ be the Laplace transform of the function $f(t)$, $$F(s)=\mathcal{L}\left(f(t)\right)(s).$$ Let $G(s)$ and $q\left(s\right)$ be two analytic functions of the complex variable $s$ such that $G(s) \exp\left(-\tau q\left(s\right)\right)$ is the Laplace transform of the function
$g\left(t, \tau\right)$,
$$G(s) \exp\left(-\tau q\left(s\right)\right)=\mathcal{L}\left(g\left(t, \tau\right)\right)(s).$$ The Efros theorem states that the following Laplace transform holds \cite{Efr35}:
$$G(s)F\left(q\left(s\right)\right) =\mathcal{L}\left(\int_0^{+\infty}
f\left( \tau\right)g\left(t, \tau\right) d \tau\right)(s).$$

The Efros theorem provides the following Laplace inversion \cite{NigMai1994}:
\begin{eqnarray}
s^{-b } F\left(s^a\right)=\mathcal{L}\left(
\int_0^{+\infty}\Phi_{a,b}\left(t,t'\right)
f\left(t'\right) dt'\right)(s).
\label {intPhiabf}
\end{eqnarray}
The function $\Phi_{a,b}\left(t,t'\right)$ be defined as follows \cite{NigMai1994}:
\begin{eqnarray}
&&\hspace{-3em}\Phi_{a,b}\left(t,t'\right)=
\pi^{-1}\int_0^t
\frac{\sin \left(t' u^a \sin\left(\pi a\right)+\pi b\right)}{u^{b}
\exp\left(u t +t' u^a \cos\left(\pi a\right)\right)}\,d u +\delta_{b,1},\label {Phiab}
\end{eqnarray}
where $\delta_{b,1}$ is the Kronecker symbol. Particularly, the function $\Phi_{a,b}\left(t,t'\right)$ reproduces the Heaviside
function, for $a=b=1$, or the delta function, for $a=1$ and $b=0$, respectively, $$\Phi_{1,1}\left(t,t'\right)= \Theta\left(t-t'\right), \hspace{1em}\Phi_{1,0}\left(t,t'\right)= \delta\left(t-t'\right).$$

The Efros theorem allows us to obtain further integral forms of the special functions under study.
 \begin{theorem}
 \label{Th3}
 The Wright function $
W_{\alpha,\beta}\left(\lambda t^{\alpha}\right)$ results to be an integral form involving the modified Bessel function of the first kind $I_0(t)$:
\begin{eqnarray}
t^{\beta-1} W_{\alpha,\beta}\left(\lambda t^{\alpha}\right)
= \int_0^{+\infty}\Phi_{\alpha,\beta-\alpha}\left(t,t_1\right)
I_0\left(2\sqrt{\lambda t_1}\right) dt_1, \label {tbm1Wab}
\end{eqnarray}
for every $t\geq 0$.
%for every $t,\alpha,\beta>0$.
\end{theorem}

{\bf Proof.}[Proof of Theorem \ref{Th3}]
Form (\ref{tbm1Wab}) is obtained by expressing the right side of Eq. (\ref{LapWright2}) as follows:
$$s^{\alpha-\beta} \frac{\exp\left(\lambda/s^{\alpha}\right)}
{s^{\alpha}}.$$
The Laplace inversion of the above-reported expression is obtained in a straightforward way from Eqs. (\ref{intPhiabf}) and (\ref{Phiab}) \cite{PrBrMa145}.
$\Box$

\begin{theorem}
\label{Th4}
The partial derivative $\partial/\partial  \alpha \left(
t^{\beta-1} W_{\alpha,\beta}\left(\lambda t^{\alpha}\right)\right)$, involving the Wright function, results to be a convoluted form with  logarithmic and other types of kernels:
\begin{eqnarray}
&& \hspace{-1em}\frac{\partial}{\partial \alpha}\,t^{\beta-1} W_{\alpha,\beta}\left(\lambda t^{\alpha}\right)
=\lambda \int_0^t  \left(\ln \left(t-t_1\right)+ \gamma_0\right)
\left(\int_0^{+\infty}\Phi_{\alpha,\beta-1}\left(t_1,t_2\right)
I_0\left(2\sqrt{\lambda t_2}\right) dt_2\right)dt_1,\nonumber \\&&
 \label{Datbm1Wab}
\end{eqnarray}
for every $t\geq 0$. %for every $t,\alpha,\beta>0$.
\end{theorem}
{\bf Proof.} [Proof of Theorem \ref{Th4}]
Form (\ref{Datbm1Wab}) is obtained by expressing the right side of Eq. (\ref{LapWright2a}) as follows:
$$-\lambda \frac{\ln s}{s}\,s^{1-\beta}\frac{\exp\left(\lambda/s^{\alpha}\right)}
{s^{\alpha}}.$$
The Laplace inversion of the above-reported expression is obtained in straightforward way from Eqs. (\ref{intPhiabf}) and (\ref{Phiab}) and the convolution theorem \cite{PrBrMa145}.
$\Box$

\begin{theorem}
\label{Th5}
The partial derivative $\partial/\partial \beta \left(
\left(t^{\beta-1} W_{\alpha,\beta}\left(\lambda t^{\alpha}\right)\right)\right)$,
involving the Wright function, results to be a convoluted form with logarithmic and other types of kernels:
\begin{eqnarray}
&& \hspace{-1.0em}
\frac{\partial}{\partial \beta}
\left(t^{\beta-1} W_{\alpha,\beta}\left(\lambda t^{\alpha}\right)\right)
=
 \int_0^t \left(\ln \left(t-t_1\right)+\gamma_0\right)
\left(\int_0^{+\infty}\Phi_{\alpha,\beta-\alpha-1}\left(t_1,t_2\right)
I_0\left(2\sqrt{\lambda t_2}\right) dt_2\right) dt_1,\nonumber \\ && \label {Dbtbm1Wab}
\end{eqnarray}
for every $t\geq 0$. %for every $t,\alpha,\beta>0$.
\end{theorem}
{\bf Proof.} [Proof of Theorem \ref{Th5}]
Form (\ref{Dbtbm1Wab}) is obtained by expressing the right side of Eq. (\ref{LapWright2b}) as follows:
$$- \frac{\ln s}{s}\,s^{\alpha-\beta+1}\frac{\exp\left(\lambda/s^{\alpha}\right)}
{s^{\alpha}}.$$
The Laplace inversion of the above-reported expression is obtained in straightforward way from Eqs. (\ref{intPhiabf}) and (\ref{Phiab}) and the convolution theorem \cite{PrBrMa145}.
$\Box$

\begin{theorem}
\label{Th6}
The partial derivative $\partial/\partial  \alpha
\left(t^{\beta-1}E_{\alpha,\beta}^{\gamma}\left(\lambda t^{\alpha}\right)\right)$, involving the Prabhakar function, results to be a convoluted form with kernel of logarithmic type:
\begin{eqnarray}
&& \hspace{-3em}
\frac{\partial}{\partial \alpha}\,
t^{\beta-1} E_{\alpha,\beta}^{\gamma}\left(\lambda t^{\alpha}\right)
=\gamma\lambda
 \int_0^t  \left(\ln \left(t-t_1\right)+ \gamma_0\right)
t_1^{\alpha+\beta-2}
E_{\alpha,\alpha+\beta-1}^{\gamma+1}\left(\lambda t_1^{\alpha}\right)
dt_1,\nonumber \\&&
 \label{Datbm1Eabg1}
\end{eqnarray}
for every $t\geq 0$. %for every $t,\alpha>0$, $\left(\alpha+\beta\right)>1$,$\gamma>-1$.
\end{theorem}

For $\alpha>0$, $\left(\alpha+\beta\right)>1$, the first term of the right side of Eq. (\ref{Datbm1Eabg1}) results to be
\begin{eqnarray}
&&\hspace{-1em}\int_0^t \left(\ln \left(t-t_1\right)\right)t_1^{\alpha+\beta-2}E^{\gamma+1}_{\alpha,\alpha+\beta-1}\left(\lambda t_1^{\alpha}\right)dt_1\nonumber \\ &&\hspace{-1em}=
t^{\alpha+\beta-1}\left(\ln t-\gamma_0\right)
\left(\left(\alpha+\beta-1\right)E^{\gamma+1}_{\alpha,\alpha+\beta-1}\left(\lambda t^{\alpha}\right)\right.
\nonumber \\ &&\hspace{-1em}\left.
+\alpha \left(\gamma+1\right) \lambda t^{ \alpha}E^{\gamma+2}_{\alpha,2\alpha+\beta-1}\left(\lambda t^{\alpha}\right)\right)\nonumber \\ &&\hspace{-1em}
-t^{\alpha+\beta-1}
\left(
\left(\alpha+\beta-1\right)
\sum_{n=0}^{\infty}\frac{\Gamma\left(n+\gamma+1\right)
\psi\left(\alpha n+\alpha+\beta\right)\left(\lambda t^{\alpha}\right)^n}{\Gamma\left(\gamma+1\right)n!\Gamma\left(\alpha n+\alpha+\beta-1\right)}\right.\nonumber \\ &&\hspace{-1em}
\left.+ \alpha \left(\gamma+1\right) \lambda t^{\alpha}
\sum_{n=0}^{\infty}\frac{\Gamma\left(n+\gamma+2\right)
\psi\left(\alpha n+2\alpha+\beta\right)\left(\lambda t^{\alpha}\right)^n}{\Gamma\left(\gamma+2\right) n!\Gamma\left(\alpha n+2\alpha+\beta-1\right)}\right)
.\label{Dbtbm1Eabg1rhs1}
\end{eqnarray}

For $\alpha>0$, $\left(\alpha+\beta\right)>1$, $\gamma>-1$, the second term of the right side of Eq. (\ref{Datbm1Eabg1}) results to be
\begin{eqnarray}
&&\hspace{-2em}\int_0^t  t_1^{\alpha+\beta-2}
E_{\alpha,\alpha+\beta-1}^{\gamma+1}\left(\lambda t_1^{\alpha}\right)
dt_1=t^{\alpha+\beta-1}E^{\gamma+1}_{\alpha,\alpha+\beta}\left(\lambda t^{\alpha}\right).
\label{Dbtbm1Eabg1rhs2}
\end{eqnarray}

{\bf Proof.} [Proof of Theorem \ref{Th6}]
Form (\ref{Datbm1Eabg1}) is obtained by expressing the right side of Eq. (\ref{Lap3parML2a}) as follows:
$$- \gamma\lambda\,\frac{\ln s}{s}\,\frac{s^{\alpha\gamma-\beta+1}
}{\left(s^{\alpha}-\lambda\right)^{\gamma+1}}.$$
The Laplace inversion of the above-reported expression is obtained by performing the Laplace inversion of each of the two fractions and applying the convolution theorem \cite{PrBrMa145}.

The term-by term integration of the series representation of the involved Prabhakar function provides Eq. (\ref{Dbtbm1Eabg1rhs1}) via Eq. (\ref{intplog}), while Eq. (\ref{Dbtbm1Eabg1rhs2}) is obtained in straightforward way via the term-by term integration \cite{PrBrMa145}.
$\Box$

\begin{theorem}
\label{Th7}
The partial derivative $\partial/\partial  \beta
\left(t^{\beta-1}E_{\alpha,\beta}^{\gamma}\left(\lambda t^{\alpha}\right)\right)$,
involving the Prabhakar function, results to be a convoluted form with kernel of logarithmic type:
\begin{eqnarray}
&& \hspace{-4em}
\frac{\partial}{\partial \beta}
\,\left(t^{\beta-1} E_{\alpha,\beta}^{\gamma}\left(\lambda t^{\alpha}\right)\right)
= \int_0^t  \left(\ln \left(t-t_1\right)+\gamma_0\right)
t_1^{\beta-2}
E_{\alpha,\beta-1}^{\gamma}\left(\lambda t_1^{\alpha}\right)
dt_1,
 \label{Dbtbm1Eabg}
\end{eqnarray}
for every $t\geq 0$. %for every $t,\alpha>0$, $\left(\alpha+\beta\right)>1$, $\gamma>-1$.
\end{theorem}

For $\alpha>0$ and $\beta>1$, the first term of the right side of Eq. (\ref{Dbtbm1Eabg}) results to be
\begin{eqnarray}
&&\hspace{-1em}\int_0^t \left(\ln \left(t-t_1\right)\right)
t_1^{\beta-2}E^{\gamma}_{\alpha,\beta-1}\left(\lambda t_1^{\alpha}\right)dt_1\nonumber \\ &&\hspace{-1em}=
t^{\beta-1}\left(\ln t-\gamma_0\right)
\left(
\left(\beta-1\right)E^{\gamma}_{\alpha,\beta-1}\left(\lambda t^{\alpha}\right)  %\right.\nonumber \\ &&\hspace{-1em}\left.
+\alpha \gamma \lambda t^{ \alpha}E^{\gamma+1}_{\alpha,\alpha+\beta-1}\left(\lambda t^{\alpha}\right)
\right)\nonumber \\ &&\hspace{-1em}
-t^{\beta-1}
\left(
\left(\beta-1\right)
\sum_{n=0}^{\infty}\frac{\Gamma\left(n+\gamma\right)
\psi\left(\alpha n+\beta\right)\left(\lambda t^{\alpha}\right)^n}{\Gamma\left(\gamma\right)n!\Gamma\left(\alpha n+\beta-1\right)}\right.\nonumber \\ &&\hspace{-1em}
\left.+ \alpha \gamma \lambda t^{\alpha}
\sum_{n=0}^{\infty}\frac{\Gamma\left(n+\gamma+1\right)
\psi\left(\alpha n+\alpha+\beta\right)\left(\lambda t^{\alpha}\right)^n}{\Gamma\left(\gamma+1\right) n!\Gamma\left(\alpha n+\alpha+\beta-1\right)}
\right)
.\label{Dbtbm1Eabgrhs1}
\end{eqnarray}

For $\alpha,\gamma>0$ and $\beta>1$, the second term of the right side of Eq. (\ref{Dbtbm1Eabg}) results to be
\begin{eqnarray}
&&\hspace{-2em}\int_0^t  t_1^{\beta-2}
E_{\alpha,\beta-1}^{\gamma}\left(\lambda t_1^{\alpha}\right)
dt_1=t^{\beta-1}E^{\gamma}_{\alpha,\beta}\left(\lambda t^{\alpha}\right).\label{Dbtbm1Eabgrhs2}
\end{eqnarray}

{\bf Proof.} [Proof of Theorem \ref{Th7}]
Form (\ref{Dbtbm1Eabg}) is obtained by expressing the right side of Eq. (\ref{Lap3parML2b}) as follows:
$$- \frac{\ln s}{s}\,\frac{s^{\alpha\gamma-\beta+1}
}{\left(s^{\alpha}-\lambda\right)^{\gamma}}.$$ The Laplace inversion of the above-reported expression is obtained by performing the Laplace inversion of each of the two fractions and applying the convolution theorem \cite{PrBrMa145}.
$\Box$

\begin{theorem}
\label{Th8}
The partial derivative $\partial/\partial \gamma \left(
t^{\beta-1}E_{\alpha,\beta}^{\gamma}\left(\lambda t^{\alpha}\right)\right)$, involving the Prabhakar function, results to be the sum of convoluted forms with logarithmic and other types of kernels:
\begin{eqnarray}
&& \hspace{-1em}
\frac{\partial}{\partial \gamma}\,
t^{\beta-1} E_{\alpha,\beta}^{\gamma}\left(\lambda t^{\alpha}\right)
= \alpha
\int_0^t\left(\ln \left(t-t_1\right)+ \gamma_0\right)
t_1^{\beta-2}E_{\alpha,\beta-1}^{\gamma}\left(\lambda t_1^{\alpha}\right)
dt_1 \nonumber \\ &&\hspace{-1em}
-\int_0^t\left(t-t_1\right)^{\beta-\alpha-1}
E_{\alpha,\beta-\alpha}^{\gamma-1}\left(\lambda \left(t-t_1\right)^{\alpha}\right) \nonumber \\ &&\hspace{-1em}\times \,
\left(\int_0^{+\infty}\Phi_{\alpha,0}\left(t_1,t_2\right)
\exp\left(\lambda t_2\right)
\left(\ln t_2+\gamma_0\right)
dt_2\right) dt_1,\nonumber \\ &&
 \label {Dgtbm1Eabg}
\end{eqnarray}
for every $t\geq 0$. %,\alpha,\gamma>0$ and $\beta>1$.
\end{theorem}

{\bf Proof.} [Proof of Theorem \ref{Th8}]
Form (\ref{Dgtbm1Eabg}) is obtained by expressing the right side of Eq. (\ref{Lap3parML2c}) as follows:
$$- \alpha  \frac{\ln s}{s}\,
\frac{s^{\alpha\gamma-\beta+1}
}{\left(s^{\alpha}-\lambda\right)^{\gamma}}
+ \frac{s^{\alpha\gamma-\beta}
}{\left(s^{\alpha}-\lambda\right)^{\gamma-1}}
\frac{\ln \left(s^{\alpha}-\lambda\right)}{\left(s^{\alpha}-\lambda\right)}
.$$
The Laplace inversion of the above-reported expression is obtained by performing the Laplace inversion of each of the two terms of the sum via the convolution theorem and Eqs. (\ref{intPhiabf}) and (\ref{Phiab}) \cite{PrBrMa145}.
$\Box$

%%%%%%%%%%%%%%%%%%%%%%%%%%%%%%%%%%%%%%%%%%
\newpage
\section{Discussion and Outlook}

The paper provides several formulas of the Laplace transforms of the derivatives with respect to parameters a number of special functions. Some of these formulas are represented via elementary or known special functions, but others are related to new types of special functions. These results could serve to improve our understanding about  the local and global behaviour of the considered special functions. They can be useful for solving differential equations.

The formulas of the Laplace transforms of the derivatives are reconstructed using the Efros theorem. Using a special case of this theorem and operational calculus, it will be possible to derive new formulas for special functions of Mittag-Leffler type and their generalizations. These formulas open new ways for considering different models, as well as for solving new classes of differential and integral equations. The obtained formulas contained series which can be considered as a new type of special functions. The coefficients in these series contain not only products of $\Gamma$-functions, but also their powers, their logarithms and their derivatives. The characteristic for these new functions is that they can be represented via the Mellin-Barnes integrals (see \cite{GKMR}). It is open the way to study asymptotics of the Laplace transforms of the derivatives of special functions with respect to parameters and to solve the corresponding differential equations.% (see \cite{GiMaRo25}).

\vs \noindent
{\bf Author contributions}
{Conceptualization, S.R. and F.M.; methodology, S.R. and F.G.;  investigation, S.R. and F.G.; data curation, S.R.; writing---original draft preparation, S.R.; writing---review and editing, F.M. All authors have read and agreed to the published version of the manuscript.}

\vs \noindent
{\bf Conflicts of interest} {The authors declare no conflicts of interest.'}

\section*{Acknowledgments}
The research by SR is partially supported by the State Program of Scientific Investigations ``Convergence-2025'', grant 1.7.01.4''.
\\
 The research activity of F.M. has been carried out in the framework of the activities of the National Group of Mathematical Physics (GNFM, INdAM). 
\\ 
 The authors are grateful to the anonymous referees for valuable suggestions which help us to improve the presentation of the results.

%%%%%%%%%%%%%%%%%%%%%%%%%%%%%%%%%%%%%%%%%%

%\reftitle{References}


% \bibliography{your_external_BibTeX_file}

%=====================================
% References, variant B: internal bibliography
%=====================================
\newpage
% ACS format
%\isAPAandChicago{}{%

\end{document}